\documentclass[12pt,a4paper,reqno]{amsart}

\usepackage{bm}
\usepackage{amssymb,amsmath,amsthm}
\usepackage{commath}
\usepackage{mathtools}
\usepackage{mathrsfs}
\usepackage{accents}
\usepackage{pgfplots}
\usepackage{accessibility}

\usepackage[utf8]{inputenc}

\usepackage{caption}
\usepackage{subcaption}

\usepackage[hyphens]{url}
\usepackage[colorlinks]{hyperref}
\hypersetup{
    allcolors = black
}
\usepackage[hyphenbreaks]{breakurl}
\usepackage[lined, linesnumbered, ruled]{algorithm2e}
\usepackage{xcolor}

\pgfplotsset{compat=1.17}

\usepackage{enumitem}
\makeatletter
\def\namedlabel#1#2{\begingroup
    #2
    \def\@currentlabel{#2}
    \phantomsection\label{#1}\endgroup
}

\newtheorem{theorem}{Theorem}[]

\newtheorem{conj}{Conjecture}

\theoremstyle{remark}
\newtheorem{remark}{Remark}

\renewcommand{\geq}{\geqslant}
\renewcommand{\leq}{\leqslant}
\renewcommand{\mod}{\, \textrm{mod }}

\newcommand{\R}{\mathbb{R}}
\newcommand{\F}{\mathbb{F}}

\newcommand{\Z}{\mathbb{Z}}

\newcommand{\eps}{\varepsilon}

\newboolean{appendix_included}
\setboolean{appendix_included}{false}

\newcommand{\appendixmacro}[1]{%
  \ifthenelse{\boolean{appendix_included}}%
    {#1}%
    {\cite[#1]{web}}%
}

\makeatletter
\@namedef{subjclassname@2020}{%
\textup{2020} Mathematics Subject Classification}
\makeatother

\author[V. Lev]{Vsevolod Lev}
\author[M. Matolcsi]{Máté Matolcsi}
\author[P.P. Pach]{P\'eter P\'al Pach}
\author[D. Varga]{Dániel Varga}

\title{On the density of Kravitz sets}

\date{\today}



\begin{document}
\baselineskip16pt

\begin{abstract}
We show that for a subset $A$ of the cyclic group of prime order $p>3$, if the sumset $A+A-2A$ is not the whole group, then $|A|\le \frac27\,p$.
Besides combinatorial arguments, we utilize a general technique involving linear programming.

\end{abstract}

\maketitle

\section{Introduction}

At the problem session of the Edinburgh meeting on additive combinatorics (2024), Noah Kravitz asked how large a subset $A\subseteq\F_p$ can be such that $A+A-2A\ne\F_p$. Here and throughout, for a prime $p$, by $\F_p$ we denote the finite field with $p$ elements, by $\lambda A$ the dilate of the set $A$ by the factor $\lambda$ 
\begin{equation*}
  \lambda A = \{ \lambda a\colon a\in A \},
\end{equation*}
and by $A+A-2A$ the sumset
\begin{equation*}
A+A-2A = \{a_1+a_2-2a_3\colon a_1,a_2, a_3\in A \}.
\end{equation*}  
(Other sumsets that occur below are defined similarly.)

Dilating and taking sumset operations can be combined to produce \emph{sums of dilates} $\lambda_1A+\dotsb+\lambda_kA$. 
In the integer settings, Bukh~\cite{Buk08} proved that for any finite set $A\subseteq \mathbb{Z}$, and any $\lambda_1,\dotsc,\lambda_k\in\Z$ with $\gcd(\lambda_1,\dotsc,\lambda_k)=1$, one has
  $$|\lambda_1A+\dots+\lambda_kA|\geq \Big(\sum\limits_{i=1}^k |\lambda_i|\Big)|A|-o(|A|),$$
where the error term $o(|A|)$ depends on the coefficients $\lambda_1,\dotsc,\lambda_k$. Pontiveros~\cite{Pon13} proved that the analogous estimate holds in $\F_p$, provided that the set $A$ has a sufficiently small density. For dense sets $A$, this may not be the case. As shown in \cite{Pon13}, for any nonzero integer $\lambda$ and any real $\eps>0$, there exists $\delta=\delta(\lambda,\eps)>0$ such that the following
holds: if $p$ is a sufficiently large prime, then there exists a set $A\subseteq  \mathbb{F}_p$ with $|A| \geq (1/2 - \eps)p$
such that $|A + \lambda A| \leq (1 - \delta)p$. Pontiveros also stated the following conjecture. 
\begin{conj}[{\cite[Conjecture~1.3]{Pon13}}]\label{j:pont}
For any nonzero integers $\lambda_1,\dots,\lambda_k \in\F_p$ and any real $\eps> 0$, there exists $\delta=\delta(\lambda,\eps)>0$ such that the following
holds: If $p$ is a sufficiently large prime, then there exists a set $A \subseteq \mathbb{F}_p$ with $|A| > (1/k -\eps)p$ such
that
$$|\lambda_1A + \dots + \lambda_kA| \leq (1 -\delta)p.$$
\end{conj}

For more results on the sums of dilates, see \cite{BS14,CHS09,CL25, HR11, HP21, Pla11}.

In this paper, we refute this conjecture. Addressing the case $\lambda_1=\lambda_2=1,\lambda_3=-2$, we prove the following result.
\begin{theorem}\label{t:main}
If $p>3$ is a prime and $A\subseteq \mathbb{F}_p$ is a subset of size $|A|>\frac27\,p$, then $A+A-2A=\mathbb{F}_p$.
\end{theorem}

Conjecturally, the assumption $|A|>\frac27\,p$ can be relaxed to $|A|>(\frac14+o(1))\,p$, while the trivial sufficient assumption is $|A|>\frac13\,p$; see a brief discussion at the end of this section.

The sumset $A+A-2A$ considered in Theorem~\ref{t:main} is particularly interesting, as it leads to the problem of finding the largest possible size of a set avoiding solutions of the equation $a_1-2a_2+a_3=c$ for some nonzero element $c\in\mathbb{F}_p$. Note that for $c=0$ one would get the problem of bounding the size of the largest 3AP-free set in $\mathbb{F}_p$.

Clearly, Theorem~\ref{t:main} refutes Conjecture~\ref{j:pont} for $k=3$, $\lambda_1=\lambda_2=1,\lambda_3=-2$. 
In fact, it is not difficult to give an explicit counterexample to Conjecture~\ref{j:pont}, and even one with all coefficients $\lambda_1,\dotsc,\lambda_k$ being pairwise distinct.
Namely, let $k=6$, $\lambda_1=10$, $\lambda_2=12$, $\lambda_3=15$, $\lambda_4=20$, $\lambda_5=30$, and $\lambda_6=-87$. We show that if $|A|>p/7>1$, then the sumset $S:=10A+12A+15A+20A+30A-87A$ is the whole group $\F_p$. Clearly, $0\in S$, and we show that $60\in 10A+12A+15A+20A+30A-87A$; applying this to the dilated set $60c^{-1}A$, we conclude that $c\in S$ for every $c\ne 0$, as well. From the assumption $|A|>p/7>1$ it follows that $(A-A)\cap \{1,2,3,4,5,6\}\ne\varnothing$; that is, the difference set $A-A$ contains at least one of the elements $1,2,3,4,5,6$. If $1\in A-A$, then we choose $a\in A$ such that also $a+1\in A$, and obtain a representation of $60$ using the identity
  $$ 10(a+1)+12a+15a+20(a+1)+30(a+1)-87a=60. $$
The cases $2\in A-A,\dotsc,6\in A-A$ can be dealt with in a similar manner, using the identities
\begin{itemize}
\item $10a+12a+15a+20a+30(a+2)-87a=60$,
\item $10a+12a+15a+20(a+3)+30a-87a=60$,
\item $10a+12a+15(a+4)+20a+30a-87a=60$,
\item $10a+12(a+5)+15a+20a+30a-87a=60$,
\item $10(a+6)+12a+15a+20a+30a-87a=60$.
\end{itemize}

\medskip

For a set $A\subset\F_p$, by $\delta(A)$ we denote the density of $A$; that is, $\delta(A)=\frac{|A|}{p}$.

For positive integers $a, b, c$ let
\begin{align*}
  M_{a,b,c}(p) 
    & := \max \{ |A|\colon A\subset \F_p, aA+bA-cA\ne \F_p \}
\intertext{and} 
   m_{a,b,c} 
     & := \limsup_{p\to \infty} \{ \delta(A)\colon  A\subset \F_p, aA+bA-cA\ne \F_p  \}.
\end{align*}
In this note we investigate the quantities $M_{1,1,2}(p)$ and $m_{1,1,2}$. We remark that the integer analogue of the latter quantity $m_{a,b,c}$
can be studied, and our methods also work in the integer settings.

If we make the identification $\F_p=\{0,1, \dots, p-1\}$, it is easy to check that the set $A:=\{1, \dotsc, \left \lfloor (p+1)/4 \right \rfloor \}$ satisfies $A+A-2A \ne \F_p$; therefore, $M_{1,1,2}(p)\ge \left \lfloor (p+1)/4 \right \rfloor $ and $m_{1,1,2}\ge 1/4$. We conjecture that this construction is the best possible for all values of $p$. We have verified this conjecture computationally up to $p\le 149$. 

Let us remark here that the trivial upper bound is $m_{1,1,2}\leq 1/3$, which can be seen as follows. If $A+A-2A\ne \mathbb{F}_p$, we may assume that $2\notin A+A-2A$ (since, if $c\ne 0$ is missing, we may replace the set $A$ by $2c^{-1}A$). Then, similarly to the above discussion, $(A-A)\cap \{1,2\}=\emptyset$ because of the representations $a+a-2a'=2$ (for $a-a'=1$) and $a+a'-2a'=2$ (for $a-a'=2$). This yields $|A|\leq p/3$ for every prime $p>2$. Note that the same bound is also a direct consequence of the Cauchy-Davenport inequality, as $|A+A-2A|\geq \min( 3|A|-2,p)$.

In the next section, we prove Theorem~\ref{t:main}
using a general linear programming approach in the spirit of~\cite{A2023}. Yet another, purely combinatorial proof is presented in Section~\ref{s:sec3}.

\section{Linear programming proof of     
   Theorem~\ref{t:main}}\label{sec2}

In this section, borrowing some ideas and notation from \cite{A2023} and \cite{M2024}, we give a linear programming proof of Theorem~\ref{t:main}. Some technical details of the proof are presented in the Appendix.

If $A\subseteq\F_p$ satisfies $A+A-2A\ne \F_p$ then, trivially, there is a nonzero element $d\in\F_p$ with $d\notin A+A-2A$. Conversely, to prove Theorem~\ref{t:main} it suffices to show that for some fixed, nonzero $d\in\F_p$, assuming $|A|>\frac27\,p$ we have $d\in A+A-2A$; applying this result to the set $c^{-1}dA$, we then conclude that $c\in A+A-2A$ for any nonzero $c\in\F_p$; that is, $A+A-2A=\F_p$. As we will see, it is convenient to choose $d=6$, but for the time being we keep the argument general. 

\medskip

In the sequel, for any set $Y\subseteq\F_p$ we use the notation $Y^1=Y$ and $Y^{-1}=\F_p\setminus Y$. 

\medskip

Following~\cite{A2023}, we will be seeking a \emph{witness set} $X=\{x_1, x_2,\dotsc,x_n\}\subset \F_p$ with certain properties. For any vector $\eps\in\{\pm1\}^n$ we associate the \emph{atomic density} 
\begin{equation}\label{atomic}
  a_X(\eps)=\delta \left ( \cap_{i=1}^n (A-x_i)^{\eps_i} \right )    
\end{equation}
and, for $I=\{i_1,\dotsc,i_r\}\subset\{1,2, \dotsc, n\}$, we define the \emph{aggregate functional} of atomic densities by      $$\overline{a}_X(I)=\sum_{\eps|_I=+1}a_X(\eps)=\delta((A-x_{i_1})\cap \dots \cap(A-x_{i_r}))$$
(where the summation extends onto all vectors $\eps=(\eps_1,\dotsc,\eps_n)$ such that $ \eps_{i_k}=1$ for each $k=1,\dots, r$). The intuition behind the atomic density $a_X(\eps)$ is the following: if we drop a randomly translated copy of $A$ in  $\F_p$, the value $a_X(\eps)$ is the probability that the intersection of this copy with $X$ is exactly the subset $\{x_j\in X: \ \eps_j=+1\}.$

Trivially, we have 
\begin{gather}
  \sum_{\eps\in \{\pm 1\}^n} a_X(\eps)=1, \label{e:C1} \\[4pt]
   a_X(\eps) \ge 0 \ \text{for all vectors}
                \ \eps\in\{\pm 1\}^n \label{e:C2},
  \intertext{and}
  \sum_{\eps\colon \eps_j=1} a_X(\eps)
    = \delta(A-x_j) = \delta(A) 
         \ \text{for all}\ 1\le j \le n. \label{e:C3}
\end{gather}
Less obviously, 
\begin{equation}\label{e:C4}
  a_X(\eps)=0 \ \text{whenever $x_i+x_j-2x_k=d$ \   
               \text{and}} \ \eps_i=\eps_j=\eps_k=1.
\end{equation}
The explanation of property \eqref{e:C4} is as follows: if $a_X(\eps)\ne 0$ while $\eps_i=\eps_j=\eps_k=1$, then the intersection $(A-x_i)\cap (A-x_j)\cap(A-x_k)$ is nonempty; consequently, there exist $a,b,c\in A$ and $y\in\F_p$ such that $y=a-x_i=b-x_j=c-x_k$ implying $x_i+x_j-2x_k=a+b-2c\ne d$. 

Similarly to~\cite{A2023, M2024}, we use translation invariance of the problem to add one more observation to our collection: if $x_{i_1},\dotsc,x_{i_r}$ is a collection of distinct elements in $X$, and $x_{j_1},\dotsc, x_{j_r}\in X$ is a translate of it, i.e. $x_{j_1}=x_{i_1}+t,\dotsc,x_{j_r}=x_{i_r}+t$ with some $t\in\F_p$, then 
\begin{equation}\label{e:C5}   
  \overline{a}_X(i_1,\dotsc,i_r)
                    =\overline{a}_X(j_1,\dotsc,j_r).
\end{equation}

We now change the viewpoint, reinterpreting $a_X(\eps)$ as formal variables of a linear program, and~\eqref{e:C1}--\eqref{e:C5} as constraints these variables satisfy. We define the objective function of our program to be $\delta(A)$, which is a linear function of the variables by equation \eqref{e:C3}; in the actual computations we use the specific choice $j=1$ in \eqref{e:C3} to define the target function. 
For each set $X\subset\F_p$, the maximum of the objective function sets an upper bound for the largest possible density $\delta(A)$. We are thus interested in finding a witness set $X=\{x_1,\dotsc,x_n\}\subset \F_p$, preferably of small cardinality, which gives a sharp upper bound on $\delta(A)$ by the linear program specified above.

We remark that, using the language of \cite{A2023, M2024}, our linear program above corresponds to the geometric fractional chromatic number (GFCN) of $X$. If the reader is familiar with those papers, this is a good analogy to keep in mind. 


\medskip

In the Appendix, we specify $d:=6$, and present a rigorous documentation of the linear program described above, for the witness set $X=\{0,2,3,4,7,8,9,10\}$. This linear program has size $22\times 53$ independently of $p$, and it proves the bound $\delta(A)\le \frac{2}{7}$ for any $p\ge 11$ (the cases, $p=5,7$ can easily be checked individually). This completes the proof of Theorem \ref{t:main}. 

It may be of some interest that, historically, we first considered the case $d=2$, the choice motivated by the argument in the Introduction, that gives $\delta(A)\le 1/3$. It turns out that for $d=2$, and for primes of the form $p=3r+1$, we can choose  $X= \{0, 1, 2, 3, r, r + 1, r + 2, 2r, 2r + 1\}$, and the resulting linear program does not depend on $r$ (which is not very hard to verify, but we will not do so here, as we will only document the case $d=6$ rigorously in the Appendix). This linear program yields the bound $\delta(A)\le 2/7$ for any $p=3r+1$.  We also observed that there exists a solution of the dual linear program with small integer coefficients. For primes of the form $p=3r-1$ we could choose $X= \{0, 1, 2, 3, r, r + 1, r + 2, 2r + 1, 2r + 3\}$, and obtain the bound $\delta(A)\le 2/7$ again. Note that the bounds above are based on the choice $d=2$, and that the witness set $X$ is concentrated around 0, $p/3$ and $2p/3$ in both cases. This suggests that if we take $d=6$,  we may be able to find a witness set $X$ concentrated around 0 only. This line of thought has led to taking $d=6$ and  $X=\{0,2,3,4,7,8,9,10\}$ as mentioned above, and as documented in the Appendix. 

In fact, whenever $n$ is co-prime to 6, this witness set $X$ testifies that for any set $A\subset \Z_n$ such that $6\notin A+A-2A$ we have $\delta(A)\le 2/7$. Thus, primality of the modulus can be replaced with the weaker assumption that the modulus is co-prime to $6$.

\begin{remark}
For small values of $p$, we can obtain better bounds by choosing the witness set $X$ appropriately.
For example, for $p=151$, $d=2$, an extensive computer search finds the $23$-element witness set shown in Figure~\ref{fig:witness_p151_n23}. It gives the upper bound $|A| \le 755/19 \approx 39.737$. Therefore, $|A|\le 39$ follows, whereas the lower bound is $38$.

\begin{figure}[ht]
    \centering
    \includegraphics[width=0.4\textwidth,trim=80 80 80 80, clip, width=.75\linewidth]{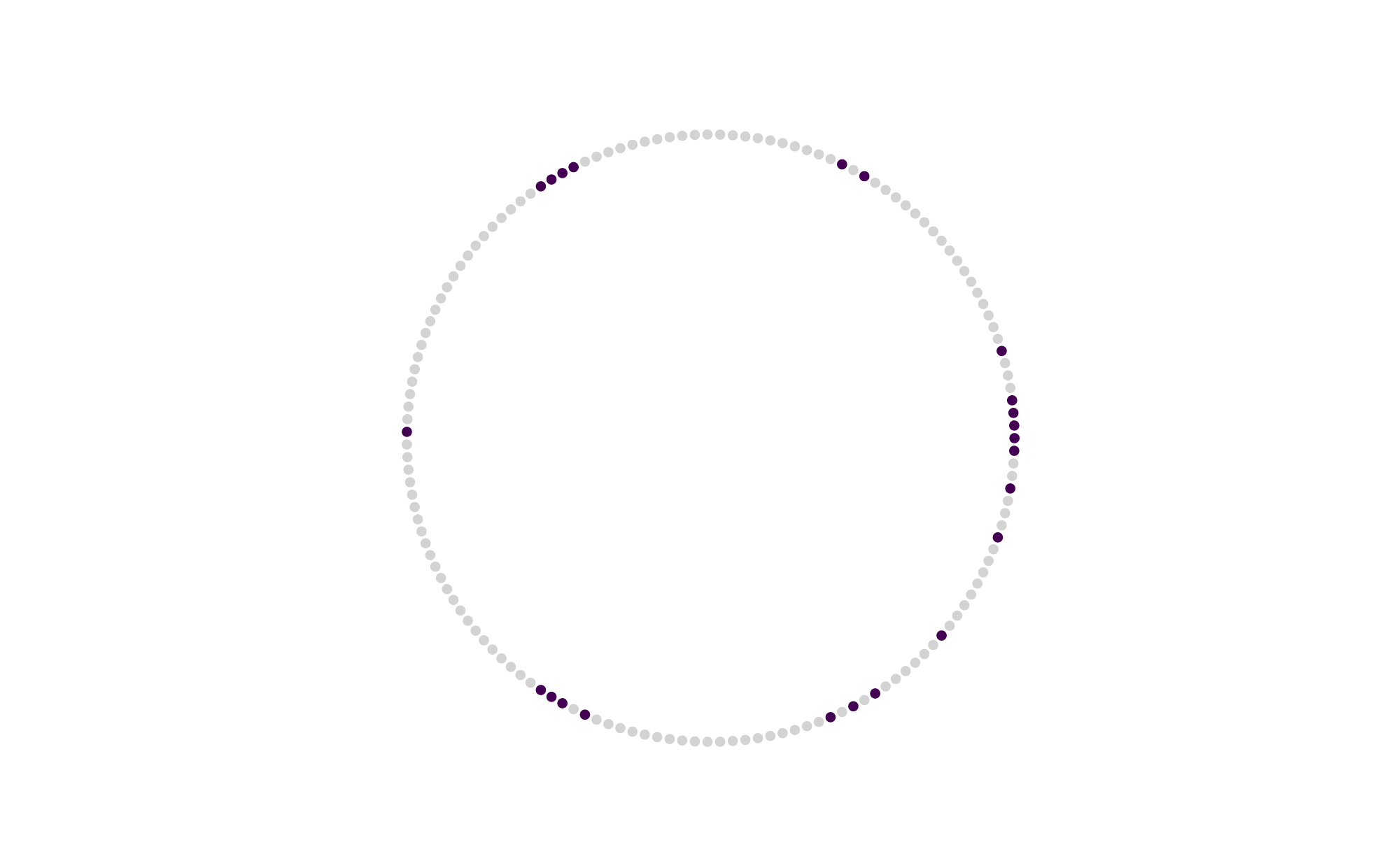}

    \caption{A 23 element witness set for $p=151$ and $d=2$, giving a density upper bound $5/19$. Removing any of the elements results in a weaker bound.}
    \label{fig:witness_p151_n23}
\end{figure}

However, it is not clear whether there is an easy way to generalize the choice of $X_p$ for all $p$. Perhaps, one can obtain the bound $\delta(A)\le 3/11$ for all $p$, by a clever choice of $X$. However, the emerging linear program testifying this bound is likely to be very large.
\end{remark}


\section{Combinatorial proof of Theorem~\ref{t:main}}\label{s:sec3}

Another way to obtain upper bounds is by means of a combinatorial argument, using ``local constraints''. 

For a nonempty subset $Y\subseteq \F_p$, let 
$$ \Delta_{d,p}(Y) := \max\{ |A\cap Y|: \ A\subseteq \F_p,\ d\notin A+A-2A \}.$$
This quantity is translation invariant: $\Delta_{d,p}(Y+t)=\Delta_{d,p}(Y)$ for any $t\in\F_p$.
Consequently, for any $A\subseteq\F_p$ with $d\notin A+A-2A$ we have 
  $$ \Delta_{d,p}(Y) \ge \frac1p\sum_{t\in\F_p}|A\cap(Y+t)| = \frac1p\,|A||Y|; $$ 
that is, 
\begin{equation}\label{e:Dfr}
    \delta(A)\le \frac{\Delta_{d,p}(Y)}{|Y|}
\end{equation}

Notice that we could have alternatively defined
$$\Delta_{d,p}(Y)= \max\{ |X|: \ X\subseteq \F_p, \ X\subseteq Y,\ d\notin X+X-2X \}, $$
and this form of the definition extends readily onto the integer settings: namely, for a set $Y\subseteq \mathbb{Z}$, let 
$$\Delta_{d,0}(Y)=\max\{ |X|: \ X\subseteq \Z, \ X\subseteq Y,\ d\notin X+X-2X \}.$$
Clearly, $\Delta_{d,0}(Y)\geq \Delta_{d,p}(Y)$ for any prime $p$, and any $Y\subseteq \F_p=\{0,1, \dots, p-1\}$.
Therefore, for any $p$, and any $A\subseteq \F_p$, equation ~\eqref{e:Dfr} implies  
  $$ \delta(A)\leq \frac{\Delta_{d,0}(Y)}{|Y|}. $$ 

In view of the linear programming argument from the previous section, we can hope to obtain a good upper bound $\frac{\Delta_{d,0}(Y)}{|Y|}$ by choosing $d=6$ and taking $Y$ to be a set concentrated around $0$. Indeed, for $d=6$ and $Y=\{0,1,\dotsc, 41\}$ one has $\Delta_{d,0}(Y)=12$. This can be verified either using a computer, or with a somewhat tedious pen-and-paper analysis, along the following lines. Suppose that $X\subset[0,41]$ is a set of size $|X|=13$ such that $6\notin X+X-2X$. Since $x-x'=3$ implies $x+x-2x'=6$, and $x-x'=6$ implies $x+x'-2x'=6$, we have $3\notin X-X$ and $6\notin X-X$. 
Therefore, each of the $14$-element sets $S_1:=\{0,3,\dots,39\}$, $S_2:=\{1,4,\dots,40\}$, and $S_3:=\{2,5,\dots,41\}$, contains at most $5$ elements of $X$. For each of these sets, there are $6$ possibilities to choose a $5$-element subset avoiding differences $3$ and $6$, and $70$, resp. $120$ possibilities to choose a $4$-element, resp. $3$-element subset satisfying this restriction. To choose $13$ elements altogether, from all three sets $S_i$, the number of chosen elements should be distributed among the sets $S_i$ as $5-5-3$ or $5-4-4$, in some order. However, all such sets turn out to contain a solution to $x+y-2z=6$ (this requires a tedious case-by-case analysis, which is possible to carry out by hand, but is much faster by computer). 

Thus, $\Delta_{6,0}(Y)=12$, implying 
  $$ \delta(A)\le \frac{\Delta_{6,0}(Y)}{|Y|} = \frac27, $$
  for any prime $p>41$.
This is the same upper bound as in Section~\ref{sec2}. For $3<p\le 41$ the bound can easily be checked for each individual value of $p$. 



An accurate case analysis also shows that most of the sets $A\subset Y$ with $|A|=12$ and $6\notin A+A-2A$ have a $\mod 7$ periodic pattern (with some irregularities due to the relatively small length 42). 

As such, in order to break the pattern mod 7, it is natural to try $d=42$ and $Y=\{0, 1, \dots, 461\}$. In this case we conjecture that $\Delta_{d,0}(Y)=126$. If so, the bound $\delta(A)\le 3/11$ would follow. Unfortunately, the case analysis becomes very long and cumbersome with such large numbers. We believe it can still be carried out with some computer assistance, and will eventually yield the bound $\delta(A)\le 3/11$. (However, we note that a simple integer programming implementation will not work: the size of the problem is far too large. One needs to do a cleverly guided case analysis.)

Again, unfortunately, we will not reach the bound 1/4 in this manner -- at least not in a computational way. As we have seen, for any prime $q$ of residue 3 mod 4, there is a construction of size $\frac{q+1}{4}$, that is of density slightly above $1/4$. This means that for any integer $m$ with $q\mid m$ the density $\frac{q+1}{4}$ can be achieved in $\mathbb{Z}_m$ (by taking a suitable mod $q$ periodic set), thus to asymptotically reach the bound $1/4$, the primeness of $p$ must be used in the argument.


\section{Concluding remarks}\label{sec4}

We proved that for any prime $p>3$, and any set $A\subseteq \mathbb{F}_p$ of size $|A|>\frac27\,p$ we have $A+A-2A=\mathbb{F}_p$. A major challenge is to improve the coefficient $\frac27$ to the conjectured $\frac14 +o(1)$. It would also be interesting to study other sums of dilates and understand how the critical density depends on the coefficients.

\medskip

Finally, we present here a Fourier analytic argument which could, in principle, lead to an improved upper bound on $\delta(A)$. The argument is a further elaboration on the linear programming approach presented in Section~\ref{sec2}, and it originates from \cite{A2023}, where it is of critical importance. However, in our present setting it only yields a very modest improvement for some values of $p$, and for some witness sets $X\subseteq \F_p$. We could not obtain any improvement independent of $p$ so far. It would be nice to understand why, and whether it can be modified to produce a significant improvement.


The function $f(x):=\delta(A\cap (A-x))= \frac{1}{p}\,1_A\ast 1_{-A}$ is positive definite; that is, can be written as 
  $$ f(x) = \frac{|A|^2}{p^2} + \sum_{m=1}^{\frac{p-1}{2}}\kappa_m\cos(2\pi mx/p), \quad x\in[0,p-1] $$
with nonnegative coefficients $\kappa_m$.  
On the other hand, if $x_i$ and $x_j$ are elements of a witness set $X=\{x_1,\dotsc,x_n\}$, with some indices $i,j\in[1,n],\ i\ne j$, then the value of $f(x_i-x_j)$ can be expressed as a sum of atomic densities:
\begin{equation*}\label{ie2}
  f(x_i-x_j) = \delta(A\cap(A-x_i+x_j)) = \overline{a}_X(\{i,j\})=\sum_{\eps_i=\eps_j=1}a_X(\eps)
\end{equation*}
where $\eps=(\eps_1,\dotsc,\eps_n)$ with $\eps_i=\eps_j=1$. Hence,
\begin{equation}\label{e:ie2}
  \frac{|A|^2}{p^2} + \sum_{m=1}^{\frac{p-1}{2}}\kappa_m\cos(2\pi m(x_i-x_j)/p) = \sum_{\eps_i=\eps_j=1}a_X(\eps). 
\end{equation}

We thus can expand our linear program with the new variables  $\kappa_m\ge 0$ (in addition to the variables $a_X(\eps)$) and new constraints \eqref{e:ie2}. Unfortunately, computations show that the resulting linear program usually has the same target value as the original one, while improvements are rare and minor.

\section{Appendix}

As claimed at the end of Section~\ref{sec2}, if $n$ is coprime to 6 and $A \subset \mathbb{Z}_n$ satisfies $6 \notin A + A - 2A$, then $\delta(A) \leq 2/7$. Here we prove this claim. All data presented here are also available online at \cite{Supplementary}, along with a computer code that verifies the required properties of the matrices defined.

The structure of Table 1 presented on the last page is as follows (the meaning of symbols $e, X, C$ etc. will be explained in the next paragraphs):

\[
\begin{array}{c|c|c}
 & & e \\ \hline
\text{congruences} & \text{dual witness } y & C
 \\
\hline
& X & \text{atomic matrix } T \\
\end{array}
\]

\medskip

The $8$-element witness set $X = \{0,2,3,4,7,8,9,10\}$ has 53 different non-trivial atomic density variables $a_X(\eps)$, i.e. variables defined by equation \eqref{atomic}, not falling into the trivial class \eqref{e:C4}. By definition, the vector $\eps$ has coordinates $\pm 1$, but in the atomic matrix $T$  below we will use the convention of changing each coordinate $-1$ to $0$. For example, the 4th column of $T$ is the vector $(0,0,0,0,0,0,1,1)$ meaning that $\eps=(-1,-1,-1,-1,-1,-1,1,1)$ in definition \eqref{atomic}. 

The vector $0\le x\in \R^{53}$ of these atomic densities satisfies the following constraints: 
$\langle \bm{1}, x \rangle = 1$ by \eqref{e:C1}, and $Cx=0$ where $C$ is the matrix corresponding to congruences defined by \eqref{e:C5}. Seemingly, we have left out the constraints \eqref{e:C3}, but note here that the constraints \eqref{e:C5} include  \eqref{e:C3} as a special case: indeed, applying \eqref{e:C5} to single element subsets of $X$ we get back \eqref{e:C3}. As an example, the topmost congruence $\{10\} \cong \{9\}$ determines the first row of the congruence matrix $C$ by taking the difference of the corresponding rows of $T$, namely the rows indexed by $10$ and $9$, respectively. Each row of $C$ is created in a similar manner: a congruence of two subsets $Y\cong Z$ is listed on the left, the rows of $T$ corresponding to the elements of $Y$ are multiplied elementwise to obtain a vector $r_{Y}$, and similarly for $Z$, and finally the difference $r_Y-r_Z$ yields the row in the matrix $C$. It is easy to see that this procedure corresponds to the equality \eqref{e:C5} in the main text. 
By convention, the sign + means a +1 entry, and -- means a $-1$ entry.

The vector $e$ is the first row of $T$ and, as such, $\langle e, x \rangle$ is the left-hand side of the constraint \eqref{e:C3} for the choice $0\in X$. 

The linear program is therefore to maximize $\langle e, x \rangle$ subject to the constraints $x\ge 0$, $\langle \bm{1}, x \rangle = 1$ and $Cx=0$. 

As such, any dual solution vector $y$ satisfying $y^T C - 7e + 2\cdot  \bm{1}\geq 0$  witnesses the fact that for any admissible vector $x$ we must have $\langle e, x \rangle \leq 2/7$. We present such a solution $y$ in the column to the left of the matrix $C$.

We note that the subsets of $X$ determine 27 congruences $Y = Z + t, Y\subseteq X, Z\subseteq X, t \in \mathbb{Z}$. Therefore, the matrix $C$ originally contained 27 rows. However, we used a greedy elimination procedure with random restarts to find a 22-element subset of the constraints that still gives the same bound $2/7$, and this is presented below for the sake of brevity.

The dual witness $y$ is not unique. The one presented here was found by a coordinate-wise iterative procedure. In each step $i$, we minimize $y_i$ with respect to $y^T C - 7 e + 2\cdot \bm{1} \geq 0$. The solution $\hat{y}_i$ turns out to be finite and an integer. We add $y_i = \hat{y}_i$ to the constraints and move on to the next coordinate $i$.
\vfill
\newpage


\[
\scalebox{0.4}{$
\begin{array}{@{}l|r|ccccccccccccccccccccccccccccccccccccccccccccccccccccc@{}}
 &  & 0 & 0 & 0 & 0 & 0 & 0 & 0 & 0 & 0 & 0 & 0 & 0 & 0 & 0 & 0 & 0 & 0 & 0 & 0 & 0 & 0 & 0 & 0 & 0 & 0 & 0 & 0 & 0 & 0 & 0 & 0 & 0 & 0 & + & + & + & + & + & + & + & + & + & + & + & + & + & + & + & + & + & + & + & + \\ \hline
\{10\} \cong \{9\} & 2 & 0 & + & - & 0 & 0 & + & - & 0 & 0 & - & 0 & - & 0 & - & 0 & - & 0 & + & 0 & + & 0 & 0 & 0 & 0 & + & - & 0 & 0 & - & 0 & - & 0 & 0 & 0 & + & - & 0 & 0 & + & - & 0 & 0 & - & 0 & - & 0 & - & 0 & - & 0 & - & 0 & - \\
\{10\} \cong \{8\} & 2 & 0 & + & 0 & + & - & 0 & - & 0 & 0 & 0 & - & - & 0 & 0 & - & - & 0 & + & - & 0 & 0 & - & 0 & 0 & + & 0 & + & 0 & 0 & 0 & 0 & 0 & 0 & 0 & + & 0 & + & - & 0 & - & 0 & 0 & 0 & - & - & 0 & 0 & - & - & 0 & 0 & 0 & 0 \\
\{10\} \cong \{7\} & 2 & 0 & + & 0 & + & 0 & + & 0 & + & - & - & - & - & 0 & 0 & 0 & 0 & 0 & + & 0 & + & - & - & 0 & 0 & + & 0 & + & - & - & 0 & 0 & 0 & 0 & 0 & + & 0 & + & 0 & + & 0 & + & - & - & - & - & 0 & 0 & 0 & 0 & 0 & 0 & - & - \\
\{10\} \cong \{4\} & 2 & 0 & + & 0 & + & 0 & + & 0 & + & 0 & 0 & 0 & 0 & - & - & - & - & 0 & + & 0 & + & 0 & 0 & - & 0 & + & 0 & + & 0 & 0 & - & - & 0 & - & 0 & + & 0 & + & 0 & + & 0 & + & 0 & 0 & 0 & 0 & - & - & - & - & 0 & 0 & 0 & 0 \\
\{10\} \cong \{3\} & -2 & 0 & + & 0 & + & 0 & + & 0 & + & 0 & 0 & 0 & 0 & 0 & 0 & 0 & 0 & - & 0 & - & 0 & - & - & - & 0 & + & 0 & + & 0 & 0 & 0 & 0 & - & - & 0 & + & 0 & + & 0 & + & 0 & + & 0 & 0 & 0 & 0 & 0 & 0 & 0 & 0 & 0 & 0 & 0 & 0 \\
\{10\} \cong \{2\} & -1 & 0 & + & 0 & + & 0 & + & 0 & + & 0 & 0 & 0 & 0 & 0 & 0 & 0 & 0 & 0 & + & 0 & + & 0 & 0 & 0 & - & 0 & - & 0 & - & - & - & - & - & - & 0 & + & 0 & + & 0 & + & 0 & + & 0 & 0 & 0 & 0 & 0 & 0 & 0 & 0 & - & - & - & - \\
\{10\} \cong \{0\} & -7 & 0 & + & 0 & + & 0 & + & 0 & + & 0 & 0 & 0 & 0 & 0 & 0 & 0 & 0 & 0 & + & 0 & + & 0 & 0 & 0 & 0 & + & 0 & + & 0 & 0 & 0 & 0 & 0 & 0 & - & 0 & - & 0 & - & 0 & - & 0 & - & - & - & - & - & - & - & - & - & - & - & - \\
\{9, 10\} \cong \{8, 9\} & -3 & 0 & 0 & 0 & + & 0 & 0 & - & 0 & 0 & 0 & 0 & - & 0 & 0 & 0 & - & 0 & 0 & 0 & 0 & 0 & 0 & 0 & 0 & 0 & 0 & + & 0 & 0 & 0 & 0 & 0 & 0 & 0 & 0 & 0 & + & 0 & 0 & - & 0 & 0 & 0 & 0 & - & 0 & 0 & 0 & - & 0 & 0 & 0 & 0 \\
\{9, 10\} \cong \{7, 8\} & -2 & 0 & 0 & 0 & + & 0 & 0 & 0 & + & 0 & 0 & - & - & 0 & 0 & 0 & 0 & 0 & 0 & 0 & 0 & 0 & - & 0 & 0 & 0 & 0 & + & 0 & 0 & 0 & 0 & 0 & 0 & 0 & 0 & 0 & + & 0 & 0 & 0 & + & 0 & 0 & - & - & 0 & 0 & 0 & 0 & 0 & 0 & 0 & 0 \\
\{9, 10\} \cong \{3, 4\} & 2 & 0 & 0 & 0 & + & 0 & 0 & 0 & + & 0 & 0 & 0 & 0 & 0 & 0 & 0 & 0 & 0 & 0 & 0 & 0 & 0 & 0 & - & 0 & 0 & 0 & + & 0 & 0 & 0 & 0 & 0 & - & 0 & 0 & 0 & + & 0 & 0 & 0 & + & 0 & 0 & 0 & 0 & 0 & 0 & 0 & 0 & 0 & 0 & 0 & 0 \\
\{9, 10\} \cong \{2, 3\} & 5 & 0 & 0 & 0 & + & 0 & 0 & 0 & + & 0 & 0 & 0 & 0 & 0 & 0 & 0 & 0 & 0 & 0 & 0 & 0 & 0 & 0 & 0 & 0 & 0 & 0 & + & 0 & 0 & 0 & 0 & - & - & 0 & 0 & 0 & + & 0 & 0 & 0 & + & 0 & 0 & 0 & 0 & 0 & 0 & 0 & 0 & 0 & 0 & 0 & 0 \\
\{8, 10\} \cong \{7, 9\} & -2 & 0 & 0 & 0 & 0 & 0 & + & 0 & + & 0 & - & 0 & - & 0 & 0 & 0 & 0 & 0 & 0 & 0 & + & 0 & 0 & 0 & 0 & 0 & 0 & 0 & 0 & - & 0 & 0 & 0 & 0 & 0 & 0 & 0 & 0 & 0 & + & 0 & + & 0 & - & 0 & - & 0 & 0 & 0 & 0 & 0 & 0 & 0 & - \\
\{8, 10\} \cong \{2, 4\} & 1 & 0 & 0 & 0 & 0 & 0 & + & 0 & + & 0 & 0 & 0 & 0 & 0 & 0 & 0 & 0 & 0 & 0 & 0 & + & 0 & 0 & 0 & 0 & 0 & 0 & 0 & 0 & 0 & - & - & 0 & - & 0 & 0 & 0 & 0 & 0 & + & 0 & + & 0 & 0 & 0 & 0 & 0 & 0 & 0 & 0 & 0 & 0 & 0 & 0 \\
\{8, 10\} \cong \{0, 2\} & 3 & 0 & 0 & 0 & 0 & 0 & + & 0 & + & 0 & 0 & 0 & 0 & 0 & 0 & 0 & 0 & 0 & 0 & 0 & + & 0 & 0 & 0 & 0 & 0 & 0 & 0 & 0 & 0 & 0 & 0 & 0 & 0 & 0 & 0 & 0 & 0 & 0 & + & 0 & + & 0 & 0 & 0 & 0 & 0 & 0 & 0 & 0 & - & - & - & - \\
\{8, 9, 10\} \cong \{7, 8, 9\} & 2 & 0 & 0 & 0 & 0 & 0 & 0 & 0 & + & 0 & 0 & 0 & - & 0 & 0 & 0 & 0 & 0 & 0 & 0 & 0 & 0 & 0 & 0 & 0 & 0 & 0 & 0 & 0 & 0 & 0 & 0 & 0 & 0 & 0 & 0 & 0 & 0 & 0 & 0 & 0 & + & 0 & 0 & 0 & - & 0 & 0 & 0 & 0 & 0 & 0 & 0 & 0 \\
\{8, 9, 10\} \cong \{2, 3, 4\} & -5 & 0 & 0 & 0 & 0 & 0 & 0 & 0 & + & 0 & 0 & 0 & 0 & 0 & 0 & 0 & 0 & 0 & 0 & 0 & 0 & 0 & 0 & 0 & 0 & 0 & 0 & 0 & 0 & 0 & 0 & 0 & 0 & - & 0 & 0 & 0 & 0 & 0 & 0 & 0 & + & 0 & 0 & 0 & 0 & 0 & 0 & 0 & 0 & 0 & 0 & 0 & 0 \\
\{4, 9\} \cong \{3, 8\} & 2 & 0 & 0 & 0 & 0 & 0 & 0 & 0 & 0 & 0 & 0 & 0 & 0 & 0 & + & 0 & + & 0 & 0 & - & - & 0 & - & 0 & 0 & 0 & 0 & 0 & 0 & 0 & 0 & + & 0 & 0 & 0 & 0 & 0 & 0 & 0 & 0 & 0 & 0 & 0 & 0 & 0 & 0 & 0 & + & 0 & + & 0 & 0 & 0 & 0 \\
\{4, 8\} \cong \{3, 7\} & 2 & 0 & 0 & 0 & 0 & 0 & 0 & 0 & 0 & 0 & 0 & 0 & 0 & 0 & 0 & + & + & 0 & 0 & 0 & 0 & - & - & 0 & 0 & 0 & 0 & 0 & 0 & 0 & 0 & 0 & 0 & 0 & 0 & 0 & 0 & 0 & 0 & 0 & 0 & 0 & 0 & 0 & 0 & 0 & 0 & 0 & + & + & 0 & 0 & 0 & 0 \\
\{4, 8, 9\} \cong \{3, 7, 8\} & -3 & 0 & 0 & 0 & 0 & 0 & 0 & 0 & 0 & 0 & 0 & 0 & 0 & 0 & 0 & 0 & + & 0 & 0 & 0 & 0 & 0 & - & 0 & 0 & 0 & 0 & 0 & 0 & 0 & 0 & 0 & 0 & 0 & 0 & 0 & 0 & 0 & 0 & 0 & 0 & 0 & 0 & 0 & 0 & 0 & 0 & 0 & 0 & + & 0 & 0 & 0 & 0 \\
\{3, 10\} \cong \{2, 9\} & -2 & 0 & 0 & 0 & 0 & 0 & 0 & 0 & 0 & 0 & 0 & 0 & 0 & 0 & 0 & 0 & 0 & 0 & + & 0 & + & 0 & 0 & 0 & 0 & 0 & - & - & 0 & - & 0 & - & 0 & 0 & 0 & 0 & 0 & 0 & 0 & 0 & 0 & 0 & 0 & 0 & 0 & 0 & 0 & 0 & 0 & 0 & 0 & - & 0 & - \\
\{3, 8, 10\} \cong \{2, 7, 9\} & 2 & 0 & 0 & 0 & 0 & 0 & 0 & 0 & 0 & 0 & 0 & 0 & 0 & 0 & 0 & 0 & 0 & 0 & 0 & 0 & + & 0 & 0 & 0 & 0 & 0 & 0 & 0 & 0 & - & 0 & 0 & 0 & 0 & 0 & 0 & 0 & 0 & 0 & 0 & 0 & 0 & 0 & 0 & 0 & 0 & 0 & 0 & 0 & 0 & 0 & 0 & 0 & - \\
\{2, 4, 9\} \cong \{0, 2, 7\} & -2 & 0 & 0 & 0 & 0 & 0 & 0 & 0 & 0 & 0 & 0 & 0 & 0 & 0 & 0 & 0 & 0 & 0 & 0 & 0 & 0 & 0 & 0 & 0 & 0 & 0 & 0 & 0 & 0 & 0 & 0 & + & 0 & 0 & 0 & 0 & 0 & 0 & 0 & 0 & 0 & 0 & 0 & 0 & 0 & 0 & 0 & 0 & 0 & 0 & 0 & 0 & - & - \\
\hline
 & 0 & 0 & 0 & 0 & 0 & 0 & 0 & 0 & 0 & 0 & 0 & 0 & 0 & 0 & 0 & 0 & 0 & 0 & 0 & 0 & 0 & 0 & 0 & 0 & 0 & 0 & 0 & 0 & 0 & 0 & 0 & 0 & 0 & 0 & 1 & 1 & 1 & 1 & 1 & 1 & 1 & 1 & 1 & 1 & 1 & 1 & 1 & 1 & 1 & 1 & 1 & 1 & 1 & 1 \\
 & 2 & 0 & 0 & 0 & 0 & 0 & 0 & 0 & 0 & 0 & 0 & 0 & 0 & 0 & 0 & 0 & 0 & 0 & 0 & 0 & 0 & 0 & 0 & 0 & 1 & 1 & 1 & 1 & 1 & 1 & 1 & 1 & 1 & 1 & 0 & 0 & 0 & 0 & 0 & 0 & 0 & 0 & 0 & 0 & 0 & 0 & 0 & 0 & 0 & 0 & 1 & 1 & 1 & 1 \\
 & 3 & 0 & 0 & 0 & 0 & 0 & 0 & 0 & 0 & 0 & 0 & 0 & 0 & 0 & 0 & 0 & 0 & 1 & 1 & 1 & 1 & 1 & 1 & 1 & 0 & 0 & 0 & 0 & 0 & 0 & 0 & 0 & 1 & 1 & 0 & 0 & 0 & 0 & 0 & 0 & 0 & 0 & 0 & 0 & 0 & 0 & 0 & 0 & 0 & 0 & 0 & 0 & 0 & 0 \\
 & 4 & 0 & 0 & 0 & 0 & 0 & 0 & 0 & 0 & 0 & 0 & 0 & 0 & 1 & 1 & 1 & 1 & 0 & 0 & 0 & 0 & 0 & 0 & 1 & 0 & 0 & 0 & 0 & 0 & 0 & 1 & 1 & 0 & 1 & 0 & 0 & 0 & 0 & 0 & 0 & 0 & 0 & 0 & 0 & 0 & 0 & 1 & 1 & 1 & 1 & 0 & 0 & 0 & 0 \\
 & 7 & 0 & 0 & 0 & 0 & 0 & 0 & 0 & 0 & 1 & 1 & 1 & 1 & 0 & 0 & 0 & 0 & 0 & 0 & 0 & 0 & 1 & 1 & 0 & 0 & 0 & 0 & 0 & 1 & 1 & 0 & 0 & 0 & 0 & 0 & 0 & 0 & 0 & 0 & 0 & 0 & 0 & 1 & 1 & 1 & 1 & 0 & 0 & 0 & 0 & 0 & 0 & 1 & 1 \\
 & 8 & 0 & 0 & 0 & 0 & 1 & 1 & 1 & 1 & 0 & 0 & 1 & 1 & 0 & 0 & 1 & 1 & 0 & 0 & 1 & 1 & 0 & 1 & 0 & 0 & 0 & 0 & 0 & 0 & 0 & 0 & 0 & 0 & 0 & 0 & 0 & 0 & 0 & 1 & 1 & 1 & 1 & 0 & 0 & 1 & 1 & 0 & 0 & 1 & 1 & 0 & 0 & 0 & 0 \\
 & 9 & 0 & 0 & 1 & 1 & 0 & 0 & 1 & 1 & 0 & 1 & 0 & 1 & 0 & 1 & 0 & 1 & 0 & 0 & 0 & 0 & 0 & 0 & 0 & 0 & 0 & 1 & 1 & 0 & 1 & 0 & 1 & 0 & 0 & 0 & 0 & 1 & 1 & 0 & 0 & 1 & 1 & 0 & 1 & 0 & 1 & 0 & 1 & 0 & 1 & 0 & 1 & 0 & 1 \\
 & 10 & 0 & 1 & 0 & 1 & 0 & 1 & 0 & 1 & 0 & 0 & 0 & 0 & 0 & 0 & 0 & 0 & 0 & 1 & 0 & 1 & 0 & 0 & 0 & 0 & 1 & 0 & 1 & 0 & 0 & 0 & 0 & 0 & 0 & 0 & 1 & 0 & 1 & 0 & 1 & 0 & 1 & 0 & 0 & 0 & 0 & 0 & 0 & 0 & 0 & 0 & 0 & 0 & 0 \\
\end{array}
$}
\]

\begin{center}
Table 1.
\end{center}

\section{Acknowledgments}

M.~M.~was supported by grants NKFIH-146387 and  and KKP 133819. 
D.~V. was supported by the Ministry of Innovation and Technology NRDI Office within the framework of the Artificial Intelligence National Laboratory (RRF-2.3.1-21-2022-00004). 
P.P.P. was supported by the Lend\"ulet program of the Hungarian Academy of Sciences (MTA) and by the National Research, Development and Innovation Office NKFIH (Grant Nr. K146387).

\newpage

\noindent
{\sc Vsevolod Lev} 

\noindent
{\em University of Haifa, Israel.}

\noindent
e-mail address: \texttt{seva@math.haifa.ac.il}

\medskip

\noindent
{\sc Máté Matolcsi} ({\em corresponding author})

\noindent
{\em HUN-REN Alfréd Rényi Institute of Mathematics, Reáltanoda u. 13-15, 1053, Budapest, Hungary,\\
and Department of Analysis and Operations Research,
Institute of Mathematics,
Budapest University of Technology and Economics,
Mûegyetem rkp. 3., H-1111 Budapest, Hungary.}

\noindent
e-mail address: \texttt{matolcsi.mate@renyi.hu}

\medskip
\noindent
{\sc P\'eter P\'al Pach}

\noindent
{\em HUN-REN Alfr\'ed R\'enyi Institute of Mathematics, Re\'altanoda utca 13--15., H-1053 Budapest,  Hungary; \\
MTA--HUN-REN RI Lend\"ulet ``Momentum'' Arithmetic Combinatorics Research Group, Re\'altanoda utca 13--15., H-1053 Budapest,  Hungary; \\
Department of Computer Science and Information Theory, Budapest University of Technology and Economics, M\H{u}egyetem rkp. 3., H-1111 Budapest, Hungary.}

\noindent
e-mail address: \texttt{pachpp@renyi.hu}

\medskip
\noindent
{\sc Dániel Varga}

\noindent
{\em HUN-REN Alfréd Rényi Institute of Mathematics, Reáltanoda u. 13-15, 1053,  Budapest, Hungary}

\noindent
e-mail address: \texttt{daniel@renyi.hu}

\end{document}